\documentclass[11pt]{amsart}

\usepackage{latexsym, graphics, psfrag, amscd, amssymb}
\newtheorem{thm}{Theorem}
\newtheorem{tm}[thm]{Theorem}
\newtheorem{fact}[thm]{Fact}
\newtheorem{lm}[thm]{Lemma}

\newtheorem{cor}[thm]{Corollary}
\newtheorem{pro}[thm]{Proposition}

\newcommand{\Si}{\Sigma}
\newcommand{\ep}{\epsilon}

\newcommand{\cM}{\mathcal{M}}

\newcommand{\Z}{\mathbb{Z}}
\newcommand{\R}{\mathbb{R}}
\newcommand{\C}{\mathbb{C}}
\newcommand{\RP}{\mathbb{RP}}

\newcommand{\Hom}{\mathrm{Hom}}

\newcommand{\fs}{\mathfrak{s}}
\newcommand{\fg}{\mathfrak{g}}
\newcommand{\ft}{\mathfrak{t}}

\newcommand{\hS}{ {\check{S}} }
\newcommand{\he}{\check{e}}

\newcommand{\tG}{ G_{ss} }
\newcommand{\tg}{\bar{g}}
\newcommand{\ta}{\bar{a}}
\newcommand{\tb}{\bar{b}}
\newcommand{\tT}{ T_{ss} }
\newcommand{\tft}{\ft_{ss} }
\newcommand{\tfg}{\fg_{ss} }

\newcommand{\bk}{\tilde{k}}

\newcommand{\mg}{{\mu^\ell_G}}

\begin{document}

\parskip=0.35\baselineskip
\baselineskip=1.2\baselineskip

\title{On the Connectedness of Moduli Spaces of Flat
Connections over Compact Surfaces}

\author{Nan-Kuo Ho}
\address{Department of Mathematics\\ National Cheng-Kung University}
\email{nankuo@mail.ncku.edu.tw}

\author{Chiu-Chu Melissa Liu}
\address{Department of Mathematics \\Harvard University}
\email{ccliu@math.harvard.edu}

\keywords{ moduli space of flat $G$ connections}

\subjclass{53}
\date{January 19, 2004}

\begin{abstract} We study the connectedness of the moduli space
of gauge equivalence classes of flat $G$-connections on a compact
orientable surface or a compact nonorientable surface for a class
of compact connected Lie groups. This class includes all the
compact, connected, simply connected Lie groups, and some
non-semisimple classical groups.
\end{abstract}
\maketitle

\section{introduction}

Given a compact Lie group $G$ and a compact surface $\Si$, let
$\cM(\Si,G)$ denote the moduli space of gauge equivalence classes
of flat $G$-connections on $\Si$. We know that $\cM(\Si,G)$ can be
identified with $\Hom(\pi_1(\Si),G)/G$, where $G$ acts on
$\Hom(\pi_1(\Si), G)$ by conjugate action of $G$ on itself (see
e.g. \cite{gold}). It is known that if $G$ is compact, connected,
simply connected (in particular, $G$ is semisimple), and $\Si$ is
orientable, then $\cM(\Si,G)$ is nonempty and connected (see e.g.
\cite{li}, \cite{gvm}, \cite{dhm}). Thus it is natural to ask about the
connectedness of $\cM(\Si,G)$ for nonorientable $\Si$. From
classification of compact surfaces, all nonorientable compact
surfaces are homeomorphic to the connected sum of the real
projective planes $\RP^2$.

Recall that we have the following structure theorem of compact connected
Lie groups \cite[Theorem 4.29]{K}:
\begin{tm}\label{thm:structure}
Let $G$ be a compact connected Lie group with center $Z(G)$,
and let $S$ be the identity component of $Z(G)$. Let
$\fg$ be the Lie algebra of $G$, and let $\tG$ be the analytic subgroup of
$G$ with Lie algebra $[\fg,\fg]$. Then $\tG$ has finite center, $S$
and $\tG$ are closed subgroups, and $G$ is the commuting product
$G=\tG S$.
\end{tm}

In Theorem \ref{thm:structure}, $G_{ss}$ is semisimple, and the
map $\tG\times S\to G=\tG S$ given by $(\tg,s)\mapsto \tg s$ is a
finite cover which is also a group homomorphism. In particular, if
$\tG$ is simply connected, we have the following result (the part
about orientable surfaces is well-known):
\begin{tm} \label{thm:torus}
Let $G$ be a compact connected Lie group with center $Z(G)$, and let
$S$ be the identity component of $Z(G)$. Let $\fg$ be the Lie algebra
of $g$, and let $\tG$ be the analytic subgroup of $G$ with Lie
algebra $[\fg,\fg]$. Suppose that $\tG$ is simply connected.
Then $\cM(\Si,G)$ is nonempty and connected if $\Si$ is a compact
orientable surface, and $\cM(\Si,G)$ is nonempty and has $2^{\dim S}$
connected components if $\Si$ is homeomorphic to $k$ copies of $\RP^2$,
where $k\neq 1,2,4$.
\end{tm}

For example, $U(n)$ and $Spin^\C(n)$ satisfy the hypothesis of Theorem
\ref{thm:torus}. We have
$$
U(n)_{ss}=SU(n),\ \
Spin^\C(n)_{ss}=Spin(n),\ \
Z(U(n))\cong Z(Spin^\C(n))\cong U(1).
$$

When $\dim S=0$, the hypothesis of Theorem \ref{thm:torus}
is equivalent to the condition that $G$ is a compact, simply connected,
connected Lie group. So we have the following special case:
\begin{cor}\label{thm:simply}
Let $G$ be a compact, connected, simply connected Lie group. Then
$\cM(\Si,G)$ is connected if $\Si$ is a compact orientable surface or
is homeomorphic to $k$ copies of $\RP^2$, where $k\neq 1,2,4$.
\end{cor}

Corollary \ref{thm:simply} is also a special case of the following result in
\cite{HL} (the part about orientable surfaces is well-known, see
\cite{li}):
\begin{thm}\label{thm:semi}
Let $G$ be a compact, connected, semisimple Lie group. Let $\Si$
be a compact orientable surface which is not homeomorphic to a
sphere. Then there is a bijection
$$
\pi_0(\cM(\Si,G))\to H^2(\Si,\pi_1(G))\cong \pi_1(G).
$$
Let $\Si$ be homeomorphic to $k$ copies of $\RP^2$, where $k\neq
1,2,4$. Then there is a bijection
$$
\pi_0(\cM(\Si,G))\to H^2(\Si,\pi_1(G))\cong \pi_1(G)/2\pi_1(G),
$$
where $2\pi_1(G)$ denote the subgroup $\{a^2\mid a\in \pi_1(G)\}$
of the finite abelian group $\pi_1(G)$.
\end{thm}

Our proofs of Theorem \ref{thm:torus}
and Theorem \ref{thm:semi} rely on the following
result of Alekseev, Malkin, and Meinrenken:
\begin{fact}[{\cite[Theorem 7.2]{gvm}}] \label{thm:fiber}
Let $G$ be a compact, connected, simply connected Lie group. Let
$\ell$ be a positive integer. Then the commutator map
$\mg:G^{2\ell}\to G$ defined by
\begin{equation} \label{eqn:commutator}
\mg(a_1,b_1,\ldots,a_\ell,b_\ell)= a_1b_1a_1^{-1}b_1^{-1}\cdots
a_\ell b_\ell a_\ell ^{-1}b_\ell^{-1}
\end{equation}
is surjective, and $(\mg)^{-1}(g)$ is connected for all $g\in G$.
\end{fact}
The surjectivity in Fact \ref{thm:fiber} follows from Goto's commutator theorem
\cite[Theorem 6.55]{HM}.

The case of orientable surfaces is discussed in Section 2. The
case of the connected sum of $2\ell+1$ copies of $\RP^2$, or
equivalently, the connected sum of a Riemann surface of genus
$\ell$ and $\RP^2$, is studied in Section 3. The case of the
connected sum of $2\ell+2$ copies of $\RP^2$, or equivalently, the
conneceted sum of a Riemann surface of genus $\ell$ and a Klein
bottle, is studied in  Section 4.

\emph{Acknowledgment}: We would like to thank Lisa Jeffrey and
Eckhard Meinrenken for many useful discussions and help.
We also wish to thank the referees for their valuable suggestions.

\section{$\Si$ is a Riemann surface with genus $\ell$} \label{sec:orientable}

In this section, we give a proof of Theorem \ref{thm:torus} for
Riemann surfaces of genus $\ell\geq 1$. The genus zero case is
trivial because the fundamental group is trivial.

\begin{lm} \label{thm:image}
Let $G, \tG, S$ be as in Theorem \ref{thm:torus}. Let $\ell$ be a
positive integer. Then the image of the commutator map
$\mg:G^{2\ell}\to G$ defined by (\ref{eqn:commutator}) is $\tG$,
and $(\mg)^{-1}(\tg)$ is connected for all $\tg\in \tG$.
\end{lm}

\paragraph{Proof.}
By Fact \ref{thm:fiber}, $\mg(\tG^{2\ell})=\tG$, so
$\mg(G^{2\ell})\supset \tG$. We now show that
$\mg(G^{2\ell})\subset \tG$. Given
$(a_1,b_1,\ldots,a_\ell,b_\ell)\in G^{2\ell}$, there exist $\ta_i,
\tb_i\in \tG$, $s_i, t_i\in S$ such that
$$
a_i=\ta_i s_i, \ \ b_i=\tb_i t_i
$$
for $i=1,\ldots,\ell$. We have
$$
\mg(a_1,b_1,\ldots,a_\ell,b_\ell)=\mg(\ta_1,\tb_1,\ldots,\ta_\ell,\tb_\ell)\in
\tG.
$$
So $\mg(G^{2\ell})\subset \tG$.

By Fact~\ref{thm:fiber}, $(\mg)^{-1}(\tg)\cap \tG^{2\ell}$ is
connected for all $\tg\in \tG$, so it suffices to show that for
any $\tg\in \tG$ and $(a_1,b_1,\ldots,a_\ell,b_\ell)\in
(\mu^\ell_G)^{-1}(\tg)$, there is a path $\gamma:[0,1]\to
(\mg)^{-1}(\tg)$ such that
$\gamma(0)=(a_1,b_1,\ldots,a_\ell,b_\ell)$ and $\gamma(1)\in
(\mg)^{-1}(\tg)\cap \tG^{2\ell}$.

Given $\tg\in \tG$ and  $(a_1,b_1,\ldots,a_\ell,b_\ell)\in
(\mg)^{-1}(\tg)$, there exist $\ta_i, \tb_i\in \tG$, $s_i, t_i\in
S$ such that
$$
a_i=\ta_i s_i, \ \ b_i=\tb_i t_i
$$
for $i=1,\ldots,\ell$. Let $\fg$ and $\fs$ be the Lie algebras of
$G$ and $S$, respectively. Then $\fs$ is a Lie subalgebra of
$\fg$. Let $\exp: \fg\to G$ be the exponential map. There exist
$X_i, Y_i\in \fs$ such that
$$
\exp(X_i)=s_i,\ \  \exp(Y_i)=t_i.
$$
for $i=1,\ldots,\ell$. Define $\gamma:[0,1]\to G^{2\ell}$ by
$$
\gamma(t)=\left(a_1\exp(-tX_1), b_1\exp(-tY_1),\ldots,
 a_l\exp(-tX_\ell), b_l\exp(-tY_\ell) \right).
$$
Then the image of $\gamma$ lies in $(\mg)^{-1}(\tg)$, and
$$
\gamma(0)=(a_1,b_1,\ldots,a_\ell,b_\ell),\ \
\gamma(1)=(\ta_1,\tb_1,\ldots,\ta_\ell,\tb_\ell)\in
(\mg)^{-1}(\tg)\cap\tG^{2\ell}. \ \ \ \Box
$$

\begin{cor}
Let $G$ be as in Theorem \ref{thm:torus}. Let $\Si$ be a Riemann
surface of genus $\ell\geq 1$. Then $\cM(\Si,G)$ is nonempty and
connected.
\end{cor}

\paragraph{Proof.} Let $\mg$ be the commutator map
defined by (\ref{eqn:commutator}),
and let $e$ be the identity element of $G$. Then
$\Hom(\pi_1(\Si),G)$ can be identified with $(\mg)^{-1}(e)$, which
is nonempty and connected by Lemma \ref{thm:image}. So
$$
\cM(\Si,G)=\Hom(\pi_1(\Si),G)/G
$$
is nonempty and connected. $\Box$

\section{$\Si$ is the connected sum of a Riemann surface of genus $\ell$
with one $\RP^2$}

The following Proposition \ref{pro:weyl} is a well-known fact. We
present an elementary proof for completeness. We use the notation
in \cite[Chapter III]{Hu}.
\begin{pro}\label{pro:weyl}
Let $\Phi$ be an irreducible root system of a Euclidean space $E$,
and let $W$ be the Weyl group of $\Phi$. Then there exists $w\in
W$ such that $1$ is not an eigenvalue of the linear map $w:E\to
E$.
\end{pro}

\paragraph{Proof.} The Coxeter element $w_c\in W$ has no eigenvalue
equal to $1$ \cite[Section 3.16]{Hu}. Here we will construct such an
element (not necessarily $w_c$) case by case.
From the classification of irreducible root systems
(see e.g. \cite[Chapter III]{Hu}), we have the following cases.

$\mathsf{A}_\ell$ ($\ell\geq 1$): $E$ is the $\ell$-dimensional
subspace of $\R^{\ell+1}$ orthogonal to the vector
$\ep_1+\cdots+\ep_{\ell+1}$, and
\[
\Phi=\{\ep_i-\ep_j \mid 1\leq i,j\leq \ell+1, i\neq j\}.
\]
The linear map $L:\R^{\ell+1}\to\R^{\ell+1}$ given by
$L\ep_i=\ep_{i+1}$ for $1\leq i\leq \ell$ and
$L\ep_{\ell+1}=\ep_1$ restricts to a linear map $w:E\to E$ which
is an element of $W$. The eigenvalues of $w$ are $e^{2\pi
i/(\ell+1)}, 1\leq i\leq \ell$.

$\mathsf{B}_\ell$ ($\ell\geq 2$): $E=\R^\ell$, and
\[
\Phi=\{\pm\ep_i \mid 1\leq i\leq \ell\}\cup
     \{\pm(\ep_i\pm\ep_j)\mid 1\leq i,j\leq \ell, i\neq j\}.
\]
The linear map $w:E\to E$ given by $v\mapsto -v$ is an element of $W$.

$\mathsf{C}_\ell$ ($\ell\geq 3$): $E=\R^\ell$, and
\[
\Phi=\{\pm 2\ep_i \mid 1\leq i\leq \ell\}\cup
     \{\pm(\ep_i\pm\ep_j)\mid 1\leq i,j\leq \ell, i\neq j\}.
\]
The linear map $w:E\to E$ given by $v\mapsto -v$ is an element of $W$.

$\mathsf{D}_\ell$ ($\ell$ is even, $\ell\geq 4$): $E=\R^\ell$, and
\[
\Phi=\{\pm(\ep_i\pm\ep_j)\mid 1\leq i,j\leq \ell, i\neq j\}.
\]
If $\ell$ is even, the linear map $w:E\to E$ given by $v\mapsto
-v$ is an element of $W$. If $l$ is odd, the linear map $w:E\to E$
given by $\ep_1\mapsto\ep_2$, $\ep_2\to -\ep_1$, and $\ep_i\mapsto
-\ep_i$ for $i\geq 3$ is an element of $W$, and the eigenvalues of
$w$ are $i,-i,-1$.

$\mathsf{E}_\ell$ ($\ell=6,8$): $E=\R^\ell$, and
\[
\Phi=\{\pm(\ep_i\pm\ep_j)\mid 1\leq i,j\leq \ell\}\cup
     \left\{\frac{1}{2}\sum_{i=1}^\ell(-1)^{k(i)}\ep_i \mid
k(i)=0,1, \sum_{i=1}^\ell k(i)\textup{ is even }\right\}.
\]
The linear map $w:E\to E$ given by $v\mapsto -v$ is an element of $W$.

$\mathsf{E}_7$: $E=\R^7$, and
\[
\Phi=\{\pm(\ep_i\pm\ep_j)\mid 1\leq i,j\leq 7\}\cup
     \left\{\frac{1}{2}\sum_{i=1}^7(-1)^{k(i)}\ep_i \mid
k(i)=0,1, \sum_{i=1}^7 k(i)\textup{ is odd }\right\}.
\]
The linear map $w:E\to E$ given by
$\ep_1\mapsto\ep_2$, $\ep_2\to -\ep_1$, and
$\ep_i\mapsto-\ep_i$ for $i\geq 3$ is an element of
$W$, and the eigenvalues of $w$ are $i,-i,-1$.

$\mathsf{F}_4$: $E=\R^4$, and
\[
\Phi=\{\pm\ep_i\mid 1\leq i\leq 4\}\cup
\{\pm(\ep_i\pm\ep_j)\mid 1\leq i,j\leq 4, i\neq j\}
\cup\left\{\pm\frac{1}{2}(\ep_1\pm\ep_2\pm\ep_3\pm\ep_4)\right\}.
\]
The linear map $w:E\to E$ given by $v\mapsto -v$ is an element of $W$.

$\mathsf{G}_2$: $E$ is the subspace of $\R^3$ orthogonal to
$\ep_1+\ep_2+\ep_3$, and
\[
\Phi=\{\pm(\ep_1-\ep_2),\pm(\ep_2-\ep_3),\pm(\ep_3-\ep_1),
 \pm(2\ep_1-\ep_2-\ep_3),\pm(2\ep_2-\ep_3-\ep_1),
 \pm(2\ep_3-\ep_1-\ep_2)\}.
\]
The linear map $L:\R^3\to\R^3$ given by $L\ep_1=\ep_2$,
$L\ep_2=\ep_3$, $L\ep_3=\ep_1$ restricts to a linear map
$w:E\to E$ which is an element of $W$.
The eigenvalues of $w$ are $e^{2\pi i/3}$ and
$e^{-2\pi i/3}$. $\Box$

The root system of a semisimple Lie algebra can be decomposed into
irreducible root systems, so Proposition \ref{pro:weyl} implies:
\begin{cor} \label{thm:weyl}
Let $G$ be a compact, connected, simply connected Lie group.
Let $\ft$ be the Lie algebra of the maximal torus $T$ of $G$.
Then there exists an element $w$ in the Weyl group of $G$
such that $1$ is not an eigenvalue of the linear map $w: \ft\to \ft$.
\end{cor}

\begin{tm} \label{thm:odd}
Let $G,\tG, S$ be as in Theorem \ref{thm:torus}. Let $\Si$ be a
compact surface whose topological type is the connected sum of a
Riemann surface of genus $\ell\geq 1$ and $\RP^2$. Then
$\cM(\Si,G)$ is nonempty and has $2^{\dim S}$ connected
components.
\end{tm}

\paragraph{Proof.}
The space $\Hom(\pi_1(\Si),G)$ can be identified with
\[
X=\{ (a_1,b_1,\ldots, a_\ell,b_\ell,c)\in G^{2\ell+1} \mid
    a_1b_1a_1^{-1}b_1^{-1}\cdots a_\ell b_\ell a_\ell ^{-1}b_\ell^{-1} c^2=e \},
\]
where $e$ is the identity element of $G$. So
$$
\cM(\Si,G) = X/G,
$$
where $G$ acts on $G^{2\ell+1}$ by diagonal conjugation. Note that
the action of $G$ preserves $X$.

Let $\tG, S$ be as in Theorem \ref{thm:torus}.
Then $\tG$ is a normal subgroup in $G$. Define
$$
\hS=G/\tG \cong S/(\tG\cap S),
$$
where the isomorphism is induced by the inclusion
$S\hookrightarrow G$. Note that $\tG\cap S\subset Z(\tG)$ is a
finite abelien group, so $S \to \hS$ is a finite cover, and $\hS$
is a compact torus with $\dim\hS=\dim S$. Let $\pi: G\to
\hS=G/\tG$ be the natural projection, and let $P:X \to G$ be
defined by
$$
(a_1,b_1,\ldots,a_\ell,b_\ell,c) \mapsto c.
$$
For $(a_1,b_1,\ldots,a_\ell,b_\ell,c)\in X$, we have
$$
c^{-2}=\mu_G(a_1,b_1,\ldots,a_\ell,b_\ell)
$$
which is an element of $\tG$ by Lemma \ref{thm:torus}.
So
$$
\pi(c)\in K \equiv\{k\in \hS\mid k^2=\he\}\cong (\Z/2\Z)^{\dim S}
$$
where $\he$ is the identity element of $\hS$. Thus $\pi\circ P$
gives a continuous map
$$
\tilde{o}: X\to K.
$$
Note that $\tilde{o}$ factors through the quotient $X/G$, so we have
$$
o:X/G \to K,
$$
and $\tilde{o}=o\circ p$, where $p$ is the projection $X\to X/G$.
Let $X_k=\tilde{o}^{-1}(k)$ for $k\in K$. We will show that $X_k$
is nonempty and connected for each $k\in K$, which implies
$o^{-1}(k)$ is nonempty and connected for each $k\in K$. This will
complete the proof since $K$ is a group of order $2^{\dim S}$.

Let $k\in K \subset \hS$. We fix $\bk\in S$ such that
$\pi(\bk)=k$. By Lemma \ref{thm:image}, $P^{-1}(c)$ is nonempty
and connected for all $c\in G$ such that $\pi(c)=k$. So $X_k$ is
nonempty. To prove that $X_k$ is connected, it suffices to show
that for any $c\in G$ such that $\pi(c)=k$, there is a path
$\gamma:[0,1]\to X$ such that $\gamma(0)\in P^{-1}(\bk)$ and
$\gamma(1)\in P^{-1}(c)$.

Since $\pi(c)=\pi(\bk)$, we have $c\bk ^{-1}\in \tG$. There exists
$g\in \tG$ such that $g^{-1}c\bk^{-1} g\in \tT$, where $\tT$ is
the maximal torus of $\tG$. Let $\tfg$ and $\tft$ denote the Lie
algebras of $\tG$ and $\tT$, respectively, and let $\exp:\tfg\to
\tG$ be the exponential map. Then
$$
g^{-1} c g \bk^{-1}=\exp(\xi)
$$
for some $\xi\in \tft$. By Corollary \ref{thm:weyl}, there exists
$w$ in the Weyl group $W$ of $\tG$ and $\xi'\in \tft$ such that
$w\cdot \xi'-\xi'=\xi$. Recall that $W=N(\tT)/\tT$, where $N(\tT)$
is the normalizer of $\tT$ in $\tG$, so $w=a\tT\in N(\tT)/\tT$ for
some $a\in\tG$. We have
$$
a\exp(t\xi') a^{-1}\exp(-t\xi')=\exp(t\xi)
$$
for any $t\in\R$.

The group $\tG$ is connected, so we may choose a path
$\tg:[0,1]\to \tG$ such that $\tg(0)=e$ and $\tg(1)=g$. Define
$\gamma:[0,1]\to G^{2\ell+1}$ by
\[
\gamma(t)=(a(t),b(t),e,\ldots,e,c(t)),
\]
where
\begin{eqnarray*}
a(t)&=&\tg(t) a\tg(t)^{-1},\\
b(t)&=&\tg(t) \exp(-2t\xi')\tg(t)^{-1},\\
c(t)&=& \bk \tg(t)\exp(t\xi)\tg(t)^{-1}.
\end{eqnarray*}
Then the image of $\gamma$ lies in $X_k$,
$\gamma(0)=(a,e,e,\ldots,e,\bk)\in P^{-1}(\bk)$, and
$$
\gamma(1)=(gag^{-1},g\exp(-2\xi')g^{-1},e,\ldots,e,c)\in P^{-1}(c).\ \ \ \Box
$$

\section{$\Si$ is the connected sum of a Riemann surface of genus $\ell\geq 2$
with a Klein bottle}

\begin{tm} \label{thm:even}
Let $G,\tG, S$ be as in Theorem \ref{thm:torus}. Let $\Si$ be a
compact surface whose topological type is the connected sum of a
Riemann surface of genus $\ell\geq 2$ and a Klein bottle. Then
$\cM(\Si,G)$ is nonempty and has $2^{\dim S}$ connected
components.
\end{tm}

\paragraph{Proof.}
The space $\Hom(\pi_1(\Si),G)$ can be identified with
\[
X=\{ (a_1,b_1,\ldots, a_\ell,b_\ell,c_1,c_2)\in G^{2\ell+2} \mid
    a_1b_1a_1^{-1}b_1^{-1}\cdots a_\ell b_\ell a_\ell ^{-1}b_\ell^{-1} c_1^2 c_2^2=e \},
\]
where $e$ is the identity element of $G$. So
$$
\cM(\Si,G) = X/G,
$$
where $G$ acts on $G^{2\ell+2}$ by diagonal conjugation. Note that
the action of $G$ preserves $X$.

Let $\pi:G\to \hS=G/\tG$ be defined as in the proof of Theorem \ref{thm:odd},
and let $P:X \to G^2$ defined by
$$
(a_1,b_1,\ldots,a_\ell,b_\ell,c_1,c_2) \mapsto (c_1,c_2).
$$
For $(a_1,b_1,\ldots,a_\ell,b_\ell,c_1,c_2)\in X$, we have
$$
c_2^{-2} c_1^{-2}=\mu_G(a_1,b_1,\ldots,a_\ell,b_\ell)
$$
which is an element of $\tG$ by Lemma \ref{thm:image}, so
$$
(\pi(c_1)\pi(c_2))^2=\pi(c_1)^2\pi(c_2)^2=\pi(c_1^2 c_2^2)=\he,
$$
where $\he$ is the the identity element of $\hS$.
Define $\tilde{o}: X\to K$ by
$$
(a_1,b_1,\ldots,a_\ell,b_\ell,c_1,c_2) \mapsto \pi(c_1)\pi(c_2)
$$
where $K$ is defined as in the proof Theorem \ref{thm:odd}.
Note that $\tilde{o}$ factors through the quotient $X/G$, so we have
$$
o:X/G \to K,
$$
and $\tilde{o}=o\circ p$, where $p$ is the projection $X\to X/G$.
Let $X_k=\tilde{o}^{-1}(k)$ for $k\in K$. We will show that
$X_k$ is nonempty and connected for each $k\in K$, which
will complete the proof.

Given $k\in K$, we fix $\bk\in S$ such that $\pi(\bk)=k$. By Lemma
\ref{thm:image}, $P^{-1}(c_1,c_2)$ is nonempty and connected for
all $(c_1,c_2)\in G^2$ such that $\pi(c_1)\pi(c_2)=k$. So $X_k$ is
nonempty. To prove that $X_k$ is connected, it suffices to show
that for any $(c_1,c_2)\in G^2$ such that $\pi(c_1)\pi(c_2)=k$,
there is a path $\gamma:[0,1]\to X$ such that $\gamma(0)\in
P^{-1}(e,\bk)$ and $\gamma(1)\in P^{-1}(c_1,c_2)$.

Given $(c_1,c_2)\in G^2$ such that $\pi(c_1)\pi(c_2)=k$, let
$k_1=\pi(c_1)$ and $k_2=\pi(c_2)$. Choose $s_1\in S$ such that
$\pi(s_1)=k_1$ and let $s_2= s_1^{-1}\bk$. Then $c_1 s_1^{-1}, c_2
s_2^{-1}\in \tG$, and  $s_1s_2=\bk$. Let $\tT,\tft,\tfg$ be as in
the proof of Theorem \ref{thm:odd}. There exist $g_1, g_2\in\tG$,
$\xi_1,\xi_2\in \tft$ such that
$$
c_1 s_1^{-1}=  g_1 \exp(\xi_1)g_1^{-1}, \ \ c_2 s_2^{-1}=  g_2
\exp(\xi_2)g_1^{-1}
$$
There exists $X\in \fs$ such that $s_1=\exp(X)$. Then
$s_2=\exp(-X)\bk$.

By Corollary \ref{thm:weyl}, there exists $w$ in the Weyl group
$W$ of $\tG$ and $\xi'_1,\xi'_2\in \tft$ such that $w\cdot
\xi_i'-\xi'_i=\xi_i$, $i=1,2$. Recall that $W=N(\tT)/\tT$, where
$N(\tT)$ is the normalizer of $\tT$ in $\tG$, so $w=a\tT\in
N(\tT)/\tT$ for some $a\in\tG$. We have
$$
a\exp(t\xi_i') a^{-1}\exp(-t\xi_i')=\exp(t\xi_i)
$$
where $i=1,2$.

The group $\tG$ is connected, so we may choose a path
$\tg_i:[0,1]\to \tG$ such that $\tg_i(0)=e$ and $\tg(1)=g_i$ for
$i=1,2$. Define $\gamma:[0,1]\to G^{2\ell+1}$ by
\[
\gamma(t)=(a_1(t),b_1(t),a_2(t),b_2(t),e,\ldots,e,c_1(t),c_2(t)),
\]
where
\begin{eqnarray*}
a_1(t)&=&\tg_2(t) a \tg_2(t)^{-1},\\
b_1(t)&=&\tg_2(t) \exp(-2t\xi'_2)\tg_2(t)^{-1}\\
a_2(t)&=&\tg_1(t) a \tg_1(t)^{-1},\\
b_2(t)&=&\tg_1(t) \exp(-2t\xi'_1)\tg_1(t)^{-1}\\
c_1(t)&=&\tg_1(t)\exp(t\xi_1)\tg_1(t)^{-1}\exp(tX)\\
c_2(t)&=&\tg_2(t)\exp(t\xi_2)\tg_2(t)^{-1}\exp(-tX)\bk \\
\end{eqnarray*}

Then the image of $\gamma$ lies in $X_k$, and
\begin{eqnarray*}
\gamma(0)&=&(a,e,a,e, e,\ldots,e, e,\bk)\in P^{-1}(e,\bk),\\
\gamma(1)&=&(g_2 a g_2^{-1},g_2\exp(-2\xi'_2)g_2^{-1},
             g_1 a g_1^{-1},g_1\exp(-2\xi'_1)g_1^{-1},
             e,\ldots, e,c_1,c_2)\in P^{-1}(c_1,c_2).\ \ \ \Box
\end{eqnarray*}

\end{document}